\documentclass{agtart_a}
\pdfoutput=1

\title[Cohomology of Coxeter groups with group ring coefficients: II]
{Cohomology of Coxeter groups\\with group ring coefficients: II} 

\author[Davis]{Michael W Davis}
\givenname{Michael W}
\surname{Davis}
\address{The Ohio State University\\
Department of Mathematics\\
231 W 18th Ave\\\newline
Columbus, Ohio 43210--1174\\
USA
}
\email{mdavis@math.ohio-state.edu}

\author[Dymara]{Jan Dymara}
\givenname{Jan}
\surname{Dymara}
\address{
Instytut Matematyczny, Uniwersytet Wroc\l awski\\ pl.
Grunwaldzki 2/4\\\newline
50-384 Wroc\l aw\\
Poland}
\email{dymara@math.uni.wroc.pl}

\author[Januszkiewicz]{Tadeusz Januszkiewicz}
\givenname{Tadeusz}
\surname{Januszkiewicz}
\address{The Ohio State University\\
Department of Mathematics\\
231 W 18th Ave\\\newline
Columbus, Ohio 43210--1174\\
USA}
\email{tjan@math.ohio-state.edu}

\author[Okun]{Boris Okun}
\givenname{Boris}
\surname{Okun}
\address{University of Wisconsin--Milwaukee\\
Department of Mathematical Sciences\\\newline
PO Box 413\\
Milwaukee WI 53201--0413\\
USA}
\email{okun@uwm.edu}

\volumenumber{6}
\issuenumber{}
\publicationyear{2006}
\papernumber{45}
\startpage{1289}
\endpage{1318}

\doi{}
\MR{}
\Zbl{}

\keyword{Coxeter group}
\keyword{Hecke algebra}
\keyword{building}
\keyword{cohomology of groups}
\subject{primary}{msc2000}{20F55}
\subject{secondary}{msc2000}{20C08}
\subject{secondary}{msc2000}{20E42}
\subject{secondary}{msc2000}{20F65}
\subject{secondary}{msc2000}{20J06}
\subject{secondary}{msc2000}{57M07}

\received{10 April 2006}
\revised{28 June 2006}
\accepted{22 June 2006}
\published{15 September 2006}
\publishedonline{15 September 2006}
\proposed{}
\seconded{}
\corresponding{}
\editor{}
\version{}

\arxivreference{math.GR/0512001}

\let\xysavmatrix\xymatrix
\def\xymatrix{\disablesubscriptcorrection\xysavmatrix}
\AtBeginDocument{\let\bar\wbar\let\tilde\wtilde\let\hat\what\def\notin{\not\in}}
\renewcommand{\theenumi}{\roman{enumi}}

\makeatletter
\def\cnewtheorem#1[#2]#3{\newtheorem{#1}{#3}[section]
\expandafter\let\csname c@#1\endcsname\c@theorem}
\makeatother

\theoremstyle{plain}
\newtheorem{theorem}{Theorem}[section]
\cnewtheorem{proposition}[theorem]{Proposition}
\cnewtheorem{conjecture}[theorem]{Conjecture}
\cnewtheorem{corollary}[theorem]{Corollary}
\cnewtheorem{lemma}[theorem]{Lemma}
\newtheorem{Theorem}{Theorem}

\newtheorem{Corollary}{Corollary}

\theoremstyle{definition}

\cnewtheorem{definition}[theorem]{Definition}
\cnewtheorem{example}[theorem]{Example}
\cnewtheorem{examples}[theorem]{Examples}
\makeautorefname{examples}{Examples}
\cnewtheorem{remark}[theorem]{Remark}
\cnewtheorem{remarks}[theorem]{Remarks}
\newtheorem{claim}{Claim}

\newtheorem{Remark}{Remark}

\newcommand{\ut}{{>T}}

\newcommand{\cac}{\mathcal {C}}

\newcommand{\cf}{\mathcal {F}}

\newcommand{\ci}{\mathcal {I}}

\newcommand{\cp}{\mathcal {P}}

\newcommand{\cs}{\mathcal {S}}

\newcommand{\cu}{\mathcal {U}}
\newcommand{\cx}{\mathcal {X}}
\newcommand{\cy}{\mathcal {Y}}

\newcommand{\qq}{{\mathbf Q}}

\newcommand{\zz}{{\mathbf Z}}

\newcommand{\hA}{{\mskip 1mu\widehat {\mskip -1mu A}}}
\newcommand{\hH}{{\widehat H}}

\newcommand{\bq}{{\mathbf q}}
\newcommand{\bs}{{\mathbf s}}

\newcommand{\wh}{\widehat}
\newcommand{\ol}{\overline}

\newcommand{\ga}{\alpha}
\newcommand{\gb}{\beta}
\newcommand{\gd}{\delta}

\newcommand{\gf}{\varphi}

\newcommand{\gO}{\Omega}

\newcommand{\gS}{\Sigma}

\newcommand{\In}{\operatorname{In}}
\newcommand{\Out}{\operatorname{Out}}
\newcommand{\Card}{\operatorname{Card}}
\newcommand{\Hom}{\operatorname{Hom}}

\newcommand{\Ker}{\operatorname{Ker}}
\newcommand{\Ima}{\operatorname{Im}}

\newcommand{\Span}{\operatorname{Span}}
\newcommand{\Res}{\operatorname{Res}}
\newcommand{\rank}{\operatorname{rank}}

\newcommand\mapright[1]{\smash{\mathop{\longrightarrow}\limits^{#1}}}

\newcommand{\comment}[1]{}

\numberwithin{equation}{section}
\begin{document}

\begin{asciiabstract}
For any Coxeter group W, we define a filtration of H^*(W;ZW) by
W-submodules and then compute the associated graded terms. More
generally, if U is a CW complex on which W acts as a reflection
group we compute the associated graded terms for H_*(U) and, in
the case where the action is proper and cocompact, for H^*_c(U).
\end{asciiabstract}

\begin{htmlabstract}
For any Coxeter group W, we define a filtration of H<sup>*</sup>(W;<b>Z</b>W) by W&ndash;submodules and then compute the associated graded
terms. More generally, if U is a CW complex on which W
acts as a reflection group we compute the associated graded terms for
H<sub>*</sub>(U) and, in the case where the action is proper and
cocompact, for H<sup>*</sup><sub>c</sub>(U).
\end{htmlabstract}

\begin{abstract}
For any Coxeter group $W$, we define a filtration of $H^*(W;{\mathbf
Z}W)$ by $W$--submodules and then compute the associated graded
terms. More generally, if $\mathcal U$ is a CW complex on which $W$
acts as a reflection group we compute the associated graded terms for
$H_*(\mathcal U)$ and, in the case where the action is proper and
cocompact, for $H^*_c(\mathcal U)$.
\end{abstract}

\maketitle

\section{Introduction}
The cohomology of a group $G$ with coefficients in a left $G$--module $M$ is denoted $H^*(G;M)$.  We are primarily interested in the case where $M$ is the group ring, $\zz G$.   Since $\zz G$ is a $G$--bimodule, $H^*(G;\zz G)$ inherits the structure of a right $G$--module.  When $G$ is discrete and acts properly and cocompactly on a contractible CW complex $\gO$,  there is a natural topological interpretation for this cohomology group:
\(
H^*(G;\zz G)\cong H^*_c(\gO)
\),
where $H^*_c(\ )$ denotes finitely supported cellular cohomology.  
The action of $G$ on $\gO$ induces a right action on cohomology and the above isomorphism is one of right $G$--modules.  For a general group $G$, not much is known about the $G$--module structure on $H^*(G;\zz G)$.  For example, even in the above case where $G$ acts properly and cocompactly on a contractible $\gO$, we don't believe it is known whether or not $H^*(G;\zz G)$ is always finitely generated as a $G$--module.

Here we deal with the case where $G=W$, a Coxeter group.  In 
\cite{d98} the first author computed $H^*(W;\zz W)$ as an abelian
group but not as a $W$--module.  We partially remedy the situation
here.  We do not quite determine the $W$--module structure on
$H^*(W;\zz W)$.  Rather, we describe a certain decreasing filtration
of $H^*(W;\zz W)$ by $W$--submodules and compute the associated graded
terms.  In order to describe this computation, we need some notation.

Suppose $(W,S)$ is a Coxeter system. 
($W$ is the group and $S$ is the distinguished set of involutions which generates $W$.)   
A  subset $T\subset S$ is \emph{spherical} if the subgroup $W_T$, generated by $T$,  is finite.  $\cs$ denotes the set of spherical subsets of $S$, partially ordered by inclusion. 

Let $A:=\zz W$ denote the group ring.
Let $\{e_w\}_{w\in W}$ be its standard basis.  For each $T\in \cs$, put 
\[
a_T:=\sum_{u\in W_T} e_u.
\]
$A^T$ denotes the right ideal $a_TA$.  If $T\subset U\in \cs$, then $a_U\in A^T$ (formula \eqref{e:aUhU} of \fullref{s:gpring}); hence, $A^U\subset A^T$.  Let $A^\ut$ be the right $W$--submodule spanned by the $A^U$, with $U\supsetneq T$. 

For each $w\in W$, put $\In'(w):=\{s\in S\mid l(sw)<l(w)\}$.  It is a fact that  $\In'(w)\in \cs$.  Set $b'_w:=a_{\In'(w)}e_w$.  We will show in \fullref{l:bbasis}  that $\{b'_w\}_{w\in W}$ is also a basis for $A$.  Define $\hA^T$ to be the $\zz$--submodule of $A$ spanned by $\{b'_w\mid w\in W, \  \In'(w)=T\}$.  N.B.\ $\hA^T$ is not a $W$--submodule of $A$; however, $\hA^T\subset A^T$ and, as we shall see in \fullref{c:hAT}, the natural map $\hA^T\to A^T/A^\ut$ is an isomorphism of free abelian groups.

The Coxeter group $W$ acts properly and cocompactly as a group generated by reflections on a certain contractible complex $\gS$ (see Davis \cite{d83}). A fundamental domain for the $W$--action on $\gS$ is a finite simplicial complex $K$, defined as  the geometric realization of the poset $\cs$.  Since $\emptyset$ is an initial element of $\cs$, $K$  is contractible.  For each $s\in S$, define $K_s$ to be the geometric realization of $\cs_{\ge\{s\}}$, where $\cs_{\ge\{s\}}:=\{T\in \cs\mid s\in T\}$.  $K_s$ is a subcomplex of $K$.  For each $U\subset S$, put $K^U:= \bigcup_{s\in U} K_s$.  The calculation of \cite{d98} was the following:
\[
H^*(W;A)=H^*_c(\gS)=\bigoplus_{T\in \cs} H^*(K,K^{S-T})\otimes \hA^T
\]
We give a new proof of  this calculation in \fullref{s:gpring} (\fullref{t:cohomU}).  This proof involves showing that a certain coefficient system on $K$ decomposes as a direct sum and that cohomology groups with coefficients in the summands are precisely the terms on the right hand side of the above formula.

A similar formula should be true for buildings.  By modifying the previous argument, we prove such a formula in \fullref{s:buildings} in the case of a locally finite, thick, right-angled building $\Phi$.  In that section, $A$ stands for the abelian group of finitely supported $\zz$--valued functions on (the set of chambers of) $\Phi$.   We define analogs of the $A^T$ and $\hA^T$ and in \fullref{t:cohombldg} we prove a formula of the form:
\[
H^i_c(|\Phi|)\cong \bigoplus_{T\in \cs} H_i(K,K^{S-T})\otimes \hA^T
\]
where $|\Phi|$ is the geometric realization of $\Phi$.

Let us  return to the question of determining the $W$--action on cohomology. 
There is a decreasing filtration of right $W$--submodules of $A$ (where $A=\zz W$):
\[
F_0\supset \cdots F_p\supset \cdots
\]
Here $F_p$ is the submodule of $A$ spanned by the $A^T$ with $\Card(T)\ge p$.  This induces a filtration of $H^*(W;A)$.  The main result of this paper is the following. (A more precise version of this result is stated as \fullref{t:main}, below.)

\begin{Theorem}
In filtration degree $p$, the graded right $W$--module associated to the above filtration of $H^*(W;A)$ is isomorphic  to
\[
\bigoplus_{\substack {T\in \cs\\ \Card (T)=p}} H^*(K,K^{S-T})\otimes (A^T/A^\ut)
\]
where $A^T/A^\ut$ is the right $W$--module defined above.
\end{Theorem}

\begin{Corollary}
$H^*(W;A)$ is finitely generated as a $W$--module.
\end{Corollary}

In view of Davis \cite{d98}, Davis et al \cite{ddjo},
Dymara and Januszkiewicz \cite{dj}, Kazhdan and Lusztig \cite{kl}
and Solomon \cite{sol}, the above computation  was a natural guess for the $W$--module structure on $H^*(W;A)=H^*_c(\gS)$.   In particular, in \cite{ddjo} we calculated the ``weighted $L^2$--cohomology'' of $\gS$ and obtained a very similar answer (provided the ``weight'' $\bq$ lies in a certain range).

We actually proceed in somewhat more generality than  indicated above.  We consider an action of $W$ as a reflection group on a CW complex $\cu$ with strict fundamental domain $X$ and then  compute certain equivariant homology and cohomology groups of $\cu$ with coefficients in $\zz W$.  The equivariant (co)homology groups we are interested in have the following well known interpretations:  $H_*^W(\cu;\zz W)=H_*(\cu)$ and when the action is proper and cocompact, $H^*_W(\cu;\zz W)=H^*_c(\cu)$.  In Theorems~\ref{t:cohomU} and \ref{t:main} we prove formulas similar to the ones above for $H^*_c(\cu)$ and $H_*(\cu)$.  In both cases the formulas involve terms of the form $H_*(X,X^U)$ or $H^*(X,X^U)$ with $U\subset S$.  The difference is that in homology only the spherical subsets $U\in \cs$ appear, while  in cohomology only  cospherical $U$ appear (ie, $S-U \in \cs$).  In the case of homology we recover the results of Davis \cite{d87}.  

When $\Phi$ is a building with a chamber transitive automorphism group $G$, one can try to calculate $H^*_c(|\Phi|)$ as a $G$--module.  We make some comments about this in the last section.
We first point out that the results of the previous sections hold when the group ring is replaced by the Hecke algebra $A_\bq$ associated to $(W,S)$ and a multiparameter $\bq$.  When $\Phi$ is a building associated to a $BN$ pair and $\bq$ its thickness, then $A_\bq$ is closely related to the algebra of finitely supported functions on $\Phi$ and hence, to the cochains on $|\Phi|$.  We  state the natural conjecture (\fullref{conj:7.4}) for the computation of $H^*_c(|\Phi|)$ as a $G$--module and prove the cochain version of it as \fullref{t:bldgcochains}. 

The first and the third authors were partially supported by NSF 
grant DMS 0405825. The second author was partially supported by KBN grant 
2 PO3A 017 25.

\section{Preliminaries}\label{s:basics}

\textbf{Invariants and coinvariants}\qua Given a left $W$--module $M$ and a subset $T\subset S$,  we have the $\zz$--submodule $M^T\subset M$  of the $W_T$--\emph{invariants} defined by
	\begin{equation}\label{e:inv}
	M^T:=M^{W_T}:=\{x\in M\mid wx=x\ \text{for all }w\in W_T\}.
	\end{equation}
More generally, for any  $\zz$--submodule $B\subset M$,  put 
	\begin{equation*}\label{e:inv2}
	B^T:=B\cap M^T. 
	\end{equation*}
	For  a right $W$--module $M$,  the $W_T$--\emph{coinvariants} are defined as a quotient $\zz$--module of $M$:
	\begin{equation}\label{e:coinv}
	M_T:=M_{W_T}:=M\otimes _{W_T} \zz \cong M/MI_T,
	\end{equation}
where $I_T$ is the augmentation ideal of $\zz W_T$ and $\zz$ is the trivial $W_T$--module.
For any  $\zz$--submodule $B\subset M$,   $B_T$ denotes the image of $B$ in $M_T$.

$\zz (W/W_T)$ denotes the (left) permutation module defined by the $W$--action  on $W/W_T$.

\begin{lemma}\label{l:invariants}
There are isomorphisms:
\begin{enumerate}
\item
$\Hom_W(\zz (W/W_T),M)\cong M^T$,
\item 
$M\otimes _W \zz(W/W_T)\cong M_T$,
\end{enumerate}
where $M$ is a left $W$--module in the first case and a right $W$--module in the second.
\end{lemma}

\begin{proof}$\phantom{99}$

(i)\qua $\Hom_W(\zz (W/W_T),M)$ can be identified with the set of $W$--equiv\-ariant functions $f\co W/W_T\to M$.  Because of equivariance, for any such $f$, $f(1W_T)\in M^T$.  Conversely, given any $x_0\in M^T$, the formula $f(wW_T)=wx_0$, gives a well-defined $f\co W/W_T\to M$.  So, $f\to f(1W_T)$ defines an isomorphism from $\Hom_W(\zz (W/W_T),M)$ to $M^T$.

(ii)\qua  We have
$$
M\otimes _W \zz(W/W_T)=M\otimes_W  \zz W\otimes _{W_T} \zz
=M\otimes _{W_T} \zz=M_T.\proved
$$
\end{proof}

\begin{remark}
If $M$ is a bimodule, then the right $W$--action on $M$ gives both\break $\Hom_W(\zz (W/W_T),M)$ and $M^T$ the structure of right $W$--modules and the isomorphism in (i) is an isomorphism of right $W$--modules.  Similarly,  (ii) is an isomorphism of left $W$--modules.  
\end{remark}

\noindent
\textbf{The basic construction}\qua Suppose $X$ is a CW complex.  Let $\cp (X)$ denote the set of cells in $X$ and  $X^{(i)}$ the set of $i$--cells. Given cells $c\in X^{(i)}$ and $c'\in X^{(i-1)}$, let $[c:c']$ denote the incidence number.  Write $c'<c$ whenever the incidence number $[c:c']$ is nonzero.  Extend this to a partial order on $\cp(X)$.

A \emph{mirror structure} on a CW complex $X$ is a family of subcomplexes
$(X_s)_{s\in S}$ indexed by some set $S$ (which for us will always be the fundamental set of generators for the Coxeter group $W$).  For each $T\subset S$, define subcomplexes of $X$:
\[
X_T:=\bigcap_{s\in T} X_s\quad\text{and}\quad X^T:=\bigcup_{s\in T} X_s 
\]
and set $X_\emptyset:=X$.  
For each  cell $c$ in $X$, set 
\[
S(c):=\{s\in S\mid c\subset X_s\}.
\]
Similarly, for each $x\in X$,  $S(x):=\{s\in S\mid x\in X_s\}$.

From the above data we construct another CW complex $\cu(W,X)$,  with a cellular $W$--action, as follows.  Give  $W$ the discrete topology. Define an equivalence relation $\sim$ on $W\times X$ by $(w,x)\sim (w',x')\iff x=x'$ and $wW_{S(x)}=w'W_{S(x')}$.  $\cu(W,X)$ is the quotient space $(W\times X)/\sim\ $.  The $W$--action on it is the obvious one.
$X$ is a fundamental domain for this action in the strict sense:  the natural inclusion $X\hookrightarrow \cu(W,X)$, which takes $x$ to the class of $(1,x)$, induces a homeomorphism $X\to \cu(W,Z)/W$.

When $X$ is the complex $K$, discussed in the Introduction, $\cu(W,K)$ is the contractible complex $\gS$.

\medskip
\textbf{Coefficient systems}\qua A \emph{system of coefficients} on a CW complex $X$ is a functor $\cf$ from $\cp(X)$ to the category of abelian groups.  Here the poset $\cp(X)$ is regarded as a category with $\Hom_{\cp(X)}(c,d)$ equal to a singleton whenever $c\leq d$ and empty otherwise. 
The functor $\cf$ will be contravariant whenever we are dealing with chains, homology or coinvariants and covariant in the case of cochains, cohomology or invariants.  Define chains and cochains with coefficients in $\cf$ by
\[
C_i(X;\cf):=\bigoplus_{c\in X^{(i)}} \cf(c)\quad\text{and}\quad 
C^i(X;\cf):=\prod_{c\in X^{(i)}} \cf(c).
\]
We regard both $i$--chains and $i$--cochains as functions $f$ from $X^{(i)}$ to $\bigcup \cf(c)$ such that $f(c)\in \cf(c)$ for each $c\in X^{(i)}$.  Boundary and coboundary maps are then defined by the usual formulas:
	\begin{align*}
	\partial (f)(c)&:=\sum [d:c] \cf_{dc} (f(d))\\
	\gd (f)(c)&:=\sum [c:d] \cf_{dc} (f(d))
	\end{align*}
where, given an $i$--cell $c$, the first sum is over all $(i+1)$--cells $d$ which are incident to $c$ and the second sum is over all $(i-1)$--cells $d$ which are incident to $c$ and where $\cf_{dc}:f(d)\to f(c)$ is the homomorphism corresponding to $d>c$ (in the first case) or $c>d$ (in the second).

\begin{examples}\label{ex:simple}
Suppose $\{X_s\}_{s\in S}$ is a mirror structure on $X$.   

(i)\qua (\textsl{Invariants})\qua Given a left $W$--module $M$, define a (covariant) system of coefficients $\ci (M)$ on $X$ by
\[
\ci(M)(c):=M^{S(c)}.
\]  
If  $B\subset M$ is any $\zz$--submodule of $M$, then we have a sub-coefficient system $\ci(B)\subset \ci(M)$, defined by $\ci(B)(c):=B^{S(c)}$.

(ii)\qua (\textsl{Coinvariants})\qua  For a right $W$--module $M$, define a  (contravariant) system of coefficients $\cac (M)$ on $X$ by 
\[
\cac(M)(c):=M_{S(c)}.
\]
Similarly, for any $\zz$--submodule $B$ of $M$, $\cac(B)(c):=B_{S(c)}$.
\end{examples}

The following observation is the key to our results.
Suppose $M$ is a left $W$--module and that we have a direct sum decomposition (of $\zz$--modules), $M=B\oplus C$, satisfying the following condition:
	\begin{equation}\label{e:Mdecomp}
	M^T=B^T\oplus C^T\quad\text{for all  $T\subset S$.}
	\end{equation}
Then we have a direct sum decomposition of coefficient systems: $\ci(M)=\ci(B)\oplus \ci(C)$.  This leads to a  decomposition of cochain groups:  $C^i(X;\ci(M))=C^i(X;\ci(B))\oplus C^i(X;\ci(C)$ and a decomposition in cohomology:		
	\begin{equation}\label{e:basic1}
	H^*(X;\ci(M))=H^*(X;\ci(B))\oplus H^*(X;\ci(C))
	\end{equation}
Similarly, suppose $M$ is a right $W$--module and $M=B\oplus C$ is a $\zz$--module decomposition satisfying:
	\begin{equation*}\label{e:Ndecomp}
	M_T=B_T\oplus C_T\quad\text{for all  $T\subset S$}
	\end{equation*}
Then we get a decomposition of coefficient systems:  $\cac(M)=\cac(B)\oplus\cac(C)$ and a corresponding decomposition of homology groups:
	\begin{equation}\label{e:basic2}
	H_*(X;\cac(M))=H_*(X;\cac(B))\oplus H_*(X;\cac(C)).
	\end{equation}

\textbf{Equivariant (co)homology}\qua Given a discrete group $G$ acting cellularly on a CW complex $\gO$, we will associate a certain equivariant homology and cohomology groups. Given a left $G$--module $M$,  the \emph{$G$--equivariant cochains on $\gO$ with coefficients in $M$} are defined by
	\begin{equation*}\label{e:C^i}
	C^i_G(\gO ;M):=\Hom_G(C_i(\gO),M).
	\end{equation*}
Similarly, if $M$ is a right $G$--module, we have the \emph{$G$--equivariant chains}
	\begin{equation*}\label{e:C_i}
	C_i^G(\gO ;M):=M\otimes_W C_i(\gO),
	\end{equation*}
where $C_i(\gO)$ denotes the group of cellular $i$--chains on $\gO$.  (Some people think that ``equivariant (co)homology''  refers to the (co)homology of $\gO\times _G EG$ with local coefficients in $M$.  However, there are  other equivariant theories, for example, the one described above.)
	
If the $G$--action is free and the projection to the orbit space is a covering projection, then equivariant (co)chains  on $\gO$ are equal to the (co)chains on the orbit space with local coefficients in $M$.  This is a useful viewpoint even when  the action is not free.  
In general, $M$ does not induce a locally constant coefficient system on the orbit space.   Rather, it induces a coefficient system on the orbit space thought of as an ``orbihedron'' or ``complex of groups.''  The theory of such coefficient systems can be found in Haefliger \cite{haefliger}.  These general coefficient systems on orbihedra are  more general then the type considered above.  (They correspond to ``lax functors'' rather than to functors.)  However, as we shall see in \fullref{l:equisheaf}, when $\gO=\cu(W,X)$, the induced coefficient system on $X$ coincides with one of the  coefficient systems described in Examples~\ref{ex:simple}. 

In the case of coefficients in the group ring,
we have the following well-known interpretation of equivariant (co)homology.

\begin{proposition}\label{p:Ugpring}
Suppose $G$ acts cellularly on a  CW complex $\gO$. Then
\begin{enumerate}
\item
$C_*^G(\gO;\zz G)\cong C_*(\gO)$.
\item
If the G-action is proper and there are only finitely many orbits of cells, then
\[
C^*_G(\gO;\zz G)\cong C^*_c(\gO).
\]
\end{enumerate}
\end{proposition}

\begin{proof}$\phantom{99}$

(i)\qua  $C_i^G(\gO;\zz G)=\zz G\otimes _{\zz G} C_i(\gO)\cong C_i(\gO)$.

(ii)\qua  For any $G$--module $M$, $\Hom_G(M,\zz G)\cong \Hom_c(M,\zz)$ (by Brown
\cite[Lemma 7.4, p age 208]{brown82}) where $\Hom_c(M,\zz)$ denotes the set of $\zz$--module maps $f\co M\to \zz$ such that for each $m\in M$, $f(gm)=0$ for all but finitely many $g\in G$.  Hence, $C^i_G(\gO;\zz G) = \Hom_G(C_i(\gO),\zz G) = \Hom_c(C_i(\gO),\zz)=C^i_c(\gO)$.
\end{proof}

Now let $\cu=\cu (W,X)$.   $W$ acts properly on  $\cu$ with compact quotient if and only if $X$ is compact and $X_U=\emptyset$ whenever $U\notin \cs$.  In view of \fullref{p:Ugpring}, \emph{when dealing with the cohomology of $\cu$, we shall always assume that these conditions hold} (ie,  $X$ is compact and $X_U=\emptyset$ for all $U\notin \cs$).  However, in the formulas for the homology of $\cu$, we need no extra assumptions on $X$.  
In the special case $\gO=\cu$, \fullref{p:Ugpring} becomes the following.

\begin{corollary}\label{cor:Ugpring}\hfil
\begin{enumerate}
\item
$C_*^W(\cu;\zz W)\cong C_*(\cu)$.
\item
$C^*_W(\cu;\zz W)\cong C^*_c(\cu)$.
\end{enumerate}
\end{corollary}

\begin{lemma}\label{l:equisheaf}\hfil
\begin{enumerate} 
\item
For any right $W$--module $M$,
$C_*^W(\cu;M)\cong C_*(X;\cac(M))$.
\item
For any left $W$--module $M$,
$C^*_W(\cu;M)\cong C^*(X;\ci(M))$.
\end{enumerate}
\end{lemma}

\begin{proof}  
Any orbit of cells in $\cu$ has the form $Wc$ for some unique cell $c$ in $X$.  As a $W$--set,  this orbit is isomorphic to $W/W_{S(c)}$.
Hence, using \fullref{l:invariants}, we get 
$$\eqalignbot{
C^i_W(\cu;M)& \cong \bigoplus_{c\in X^{(i)}} M^{S(c)}\ =C^i(X;\ci(M))\cr
C_i^W(\cu;M)&\cong \bigoplus_{c\in X^{(i)}} M_{S(c)}\ =C_i(X;\cac(M)).}
\proved$$
\end{proof}

\begin{Remark} 
The isomorphisms in \fullref{cor:Ugpring}~(ii) and \fullref{l:equisheaf}~(ii) give
\[
C^i_c(\cu)\cong C^i_W(\cu;\zz W)\cong C^i(X;\ci(\zz W)) \cong
\bigoplus_{c\in X^{(i)}} (\zz W)^{S(c)}.
\]
The composition of these gives an isomorphism $C^i_c(\cu)\to \bigoplus_{c\in X^{(i)}} (\zz W)^{S(c)}$, such that its component corresponding to $c\in X^{(i)}$ is given by
\[
f \to \sum_{w\in W} f(w^{-1}c)e_w
\]
where $(e_w)$ is the standard basis for $\zz W$.  Similarly, the composition of the  isomorphisms in \fullref{cor:Ugpring}~(i) and \fullref{l:equisheaf}~(i)  gives the obvious identification 
\[
C_i(\cu)\cong \bigoplus_{c\in X^{(i)}} (\zz W)_{S(c)}.
\]
\end{Remark}

\section{Group ring coefficients}\label{s:gpring}
\textbf{Subsets of $\boldsymbol{W}$}\qua
For any  $U\subset S$, put
\begin{align*}
\cx_U&:=\{w\in W\mid l(sw)>l(w) \text{ for all } s\in U\}\\ \cy_U&:=\{w\in W\mid l(ws)>l(w) \text{ for all } s\in U\}=(\cx_U)^{-1}
\end{align*}
$\cx_U$ (resp.\ $\cy_U$) is the set of elements in $W$ which are \emph{$(U,\emptyset)$--reduced} (resp.\ \emph{$(\emptyset,U)$--reduced}).  $\cx_U$ (resp.\ $\cy_U$) is a set of representatives for $W_U\backslash W$ (resp.\ $W/W_U$).

Given $w\in W$, set 
	\begin{align*}
	\In (w)&:=\{s\in S\mid l(ws)< l(w)\},\\ 
	\In' (w)&:=\{s\in S\mid l(sw)< l(w)\}=\In(w^{-1}). 
	\end{align*}
$\In(w)$ (resp.\ $\In'(w)$) is the set of letters of $S$ with which a reduced word for $w$ can end (resp.\ begin).
By \cite[Lemma 7.12]{d83}, for any $w\in W$, $\In (w)$ is a spherical subset.  We note that, for any $T\in\cs$,
	\begin{align*}
	w_T\cx_T&=\{w\in W\mid T\subset \In'(w)\},\\ 
	\cy_{S-T}&=\{w\in W\mid T\supset \In(w)\}, 
	\end{align*}
where $w_T\in W_T$ is  the element of longest length.  Thus, $w_T\cx_T$ is also a set of representatives for $W_T\backslash W$.

\medskip
\textbf{Symmetrization and alternation}\qua  From now on, except in \fullref{s:buildings}, $A$ denotes the group ring $\zz W$.  
For each spherical subset $T$ of $S$, define elements $a_T$ and $h_T$ in $A$ by
	\begin{equation}\label{e:aThT}
	a_T:=\sum_{w\in W_T} e_w\quad \text{and} \quad h_T:=\sum_{w\in W_T} (-1)^{l(w)}e_w
	\end{equation}
called, respectively, \emph{symmetrization} and \emph{alternation} with respect to $T$.  If $T\subset U\in \cs$, define
\[
c_{(U,T)}:=\sum_{u\in \cx_T\cap W_U} e_u\quad\text{and}\quad
d_{(U,T)}:=\sum_{u\in \cy_T\cap W_U} (-1)^{l(w)}e_u.
\]
It is easily checked that 
	\begin{equation}\label{e:aUhU}
	a_U=a_Tc_{(U,T)}\quad\text{and}\quad 
	h_U=d_{(U,T)}h_T.
	\end{equation}
For any subset $T$ of $S$, let $A^T$ denote the $W_T$--invariants in $A$, defined as in \eqref{e:inv}.  Notice that $A^T$ is $0$ if $T\notin \cs$ and is equal to the right ideal $a_TA$ if $T\in \cs$.  Similarly, for $T\in \cs$, define  $H^T$ to be the left ideal $Ah_T$ and to be $0$ otherwise.  By \eqref{e:aUhU}, $A^U\subset A^T$ and $H^U\subset H^T$ whenever $T\subset U$.  
Let $A_T$ denote the $W_T$--coinvariants, defined as in \eqref{e:coinv} and  let $I_T$ denote the augmentation ideal of $\zz W_T$. For any $s\in S$,  note that  $AI_{\{s\}}=H^{\{s\}}$. Hence, $A_{\{s\}}=A/H^{\{s\}}$.  More generally, for any $T\subset S$,
	\begin{equation*}\label{e:augment}
	AI_T=\sum_{s\in T}H^{\{s\}}\quad \text{so,}\quad
	A_T=A/\sum_{s\in T}H^{\{s\}}.
	\end{equation*}

\noindent
\textbf{Two bases for $\boldsymbol{A}$}\qua
For each $w\in W$, define elements  $b'_w, b_w\in A$ by
\begin{equation*}\label{e:bw}
b'_w:= a_{\In'(w)} e_w\quad\text{and}\quad b_w:=e_w h_{\In(w)}.
\end{equation*}

\begin{lemma}\label{l:bbasis}
$\{b'_w\mid w\in W\}$ is a basis for $A$ (as a  $\zz$--module).  More generally, for any $T\in\cs$, $\{b'_w\mid T\subset \In'(w)\}$ is a basis for $A^T$.
\end{lemma}

\begin{proof}
We first show $\{b'_w\mid w\in W\}$ is a basis.  The point is that the matrix which expresses the $b'_w$ in terms of the $e_w$ has $1$'s on the diagonal and is ``upper triangular with respect to word length.''  In detail:  first note that $b'_v$ is the sum of $e_v$ with various $e_w$ having $l(w)<l(v)$. 
Suppose $\sum \gb_w b'_w=0$ is a nontrivial linear relation.  Let $v\in W$ be an element with $\gb_v\neq 0$ and $l(v)$ maximum.  Since the coefficient of $e_v$ in the linear relation must be $0$, we have $\gb_v=0$, a contradiction.  Similarly, one shows, by induction on word length, that each $e_v$ is a linear combination of $b'_w$ with $l(w)\leq l(v)$.  Hence, $\{b'_w\}$ spans $A$.

To prove the second sentence, we must first show that $b'_w\in A^T$ whenever $T\subset \In'(w)$.  If this condition holds, then, by \eqref{e:aUhU}, 
\[
b'_w=a_{\In'(w)}e_w=a_Tc_{(\In'(w),T)}e_w \in A^T.
\] 
Note that $T\subset \In'(w)$ if and only if $w\in w_T\cx_T$.
Since, by the previous paragraph, $\{b'_w\mid w\in w_T\cx_T\}$ is linearly independent, it remains to show that it spans $A^T$. Since $w_T\cx_T$ is a set of coset representatives for $W_T\backslash W$,
a basis for $A^T$ is $\{a_Te_w\mid w\in w_T\cx_T\}$.  Let $\ol{e}_w:=c_{(\In'(w),T)}e_w$.  For $w\in w_T\cx_T$,  the matrix which expresses $\{\ol{e}_w\mid w\in w_T\cx_T\}$ in terms of $\{e_w\mid w\in w_T\cx_T\}$ has $1$'s on the diagonal and is upper triangular with respect to word length.  So,
\[
\{a_T\ol{e}_w\mid w\in w_T\cx_T\}=\{b'_w\mid T\subset \In'(w)\}
\]
is also a basis for $A^T$.
\end{proof}

\begin{lemma}\label{l:b'basis}
$\{b_w\mid w\in W\}$ is a basis for $A$.  More generally, for any subset $U$ of $S$, the projection $A\to A_{S-U}$  maps $\{b_w\mid U\supset \In(w)\}$ injectively to a basis for $A_{S-U}$.
\end{lemma}

\begin{proof}
The proof of the first sentence is omitted since it is similar to that of the first sentence of the previous lemma.

Fix a  subset $U\subset S$ and let $p\co A\to A_{S-U}$ denote the projection.  Since $A_{S-U}=\zz (W/W_{S-U})$, $\{p(e_w)\mid w\in \cy_{S-U}\}$ is the  obvious basis for $A_{S-U}$ (as a $\zz$--module).   Any element $y\in A$ can be written in the form
	\begin{equation}\label{e:y}
	y=\sum_{w\in \cy_{S-U}}\ \sum_{u\in W_{S-U}} \ga_{wu} e_{wu}.
	\end{equation}
Moreover, $y\in AI_{S-U}=\Ker (p)$ if and only if $\sum_{u\in W_{S-U}} \ga_{wu}=0$ for each $w\in \cy_{S-U}$.  Let $y$ be an element in  the submodule spanned by $\{b_w\mid U\supset \In(w)\}$ ($=\{b_w\mid w\in \cy_{S-U}\}$), ie, let
\[
y=\sum_{w\in \cy_{S-U}}y_w b_w.
\]
Suppose $p(y)=0$.  Let $v\in \cy_{S-U}$ be such that $y_v\neq 0$ and $l(v)$ is maximum with respect to this property.  Since $b_v$ is the sum of $e_v$ and $\pm 1$ times various $e_w$ with $l(w)<l(v)$, the coefficients $\ga_{vu}$  in \eqref{e:y} are $0$ for all $u\neq 1$ in $W_{S-U}$.  So, $\sum \ga_{vu} =0$ forces $\ga_{vu}=0$, a contradiction.  Thus, $\{p(b_w)\mid w\in \cy_{S-U}\}$ is linearly independent in $A_{S-U}$.    The usual argument, using induction on word length,  shows that $\{p(b_w)\mid w\in Y_{S-U}\}$ spans $A_{S-U}$.
\end{proof}

In view of Lemmas~\ref{l:bbasis} and \ref{l:b'basis},  for each $T\in \cs$, we define $\zz$--submodules of $A$:
	\begin{align*}
	\hA^T:=\Span \{b'_w\mid \In'(w)=T\},\\
	\hH^T:=\Span \{b_w\mid \In(w)=T\}.
	\end{align*}
A corollary to \fullref{l:bbasis} is the following.
\begin{corollary}\label{c:hAT}
For any $U\in\cs$,
	\[
	A^U= \bigoplus_{T\in \cs_{\geq U}} \hA^T .
	\]
\end{corollary}
Consequently, given $T\in \cs$, for any $U\subset S$ we have:
	\begin{equation}\label{e:hAT1}
	(\hA^T)^U=
	\begin{cases}
	\hA^T,	&\text{if $U\subset T$;}\\
	0,	&\text{if $U\cap (S-T)\neq \emptyset$.}
	\end{cases}
	\end{equation}
It follows that the direct sum decomposition in \fullref{c:hAT} satisfies \eqref{e:Mdecomp} and hence, gives a decomposition of coefficient systems:
	\begin{equation}\label{e:Asheaf}
	\ci(A)=\bigoplus_{T\in\cs} \ci(\hA^T).
	\end{equation}
In terms of the $\hH^T$, the version of this we are interested in is
the following:
\begin{equation}\label{e:hHT1}
	(\hH^T)_U\cong
	\begin{cases}
	\hH^T,	&\text{if $U\subset S-T$;}\\
	0,	&\text{if $U\cap T\neq \emptyset$.}
	\end{cases}
	\end{equation}
In the above  formula, by writing $(\hH^T)_U\cong \hH^T$, we mean that the projection $A\to A_U$ maps  $\hH^T$ isomorphically onto $(\hH^T)_U$.  To see that $(\hH^T)_U=0$ when $U\cap T\neq \emptyset$, note that if $s\in T\cap U$, then $\hH^T\subset H^s\subset AI_U$.

The $\hH^T$ version of \fullref{c:hAT} is the following.

\begin{corollary}\label{c:hHT}
For any $U\subset S$,
\[
A_{S-U}= \bigoplus_{T\in \cs_{\leq U}} (\hH^T)_{S-U} .
\]
\end{corollary}

So, the decomposition in \fullref{c:hHT} gives a decomposition of coefficient systems:
	\begin{equation}\label{e:Acosheaf}
	\cac(A)=\bigoplus_{T\in\cs} \cac(\hH^T).
	\end{equation}
Hence, \eqref{e:basic1} and \eqref{e:basic2} apply to give the following calculation of (co)homology with group ring coefficients.

\begin{theorem}\label{t:cohomU} Let $\cu=\cu(W,X)$. Then
\begin{align*}
H^i_c(\cu)&\cong \bigoplus_{T\in \cs} H^i(X,X^{S-T})\otimes \hA^T,\\ 
H_i(\cu)&\cong \bigoplus_{T\in \cs} H_i(X,X^T)\otimes \hH^T.
\end{align*}
\end{theorem}

\begin{proof}
To prove the first formula, 
note that by \fullref{p:Ugpring} and observation \eqref{e:basic1}, 
\[
C^i_c(\cu)=C^i(X;\ci(A))=\bigoplus_{T\in\cs}C^i(X;\ci(\hA^T)).
\]
Given a cell $c\in X^{(i)}$, by \eqref{e:hAT1}, 
	\[
	(\hA^T)^{S(c)}=
	\begin{cases}
	0,	&\text{if $c\subset X^{S-T}$;}\\
	\hA^T, &\text{otherwise.}
	\end{cases}
	\]
Hence, 
\[
C^i(X;\ci(\hA^T))=\{f\co X^{(i)}\to \hA^T\mid f(c)=0 \text{ if } c\subset X^{S-T}\}
=C^i(X,X^{S-T})\otimes \hA^T.
\]
Combining these formulas and taking cohomology, we get the first formula.

To prove the second formula,  
note that by \fullref{p:Ugpring} and observation \eqref{e:basic2}, 
\[
C_i(\cu)=C_i(X;\cac(A))=\bigoplus_{T\in\cs}C_i(X;\cac(\hH^T)).
\]
Given a cell $c\in X^{(i)}$, by \eqref{e:hHT1}, 
	\[
	(\hH^T)_{S(c)}\cong
	\begin{cases}
	0,	&\text{if $c\subset X^T$;}\\
	\hH^T, &\text{otherwise.}
	\end{cases}
	\]
Hence,
\[
C_i(X;\cac(\hH^T))=\bigoplus_{\substack{c\in X^{(i)}\\c\not\subset X^T} }\hH^T\cong C_i(X,X^T)\otimes \hH^T.
\]
Taking homology, we get the second formula. 
\end{proof}

\begin{Remark}
The first formula in \fullref{t:cohomU} is one of the main results of \cite{d98}.  (Actually, only a special case is stated in \cite{d98}; however, the general result is stated in \cite{dmei}.)  The second formula is  the main result of \cite{d87}.
\end{Remark}

\section{The $W$--module structure of $H^*_c(\cu)$ and $H_*(\cu)$}\label{s:w}

$A$ is a $W$--bimodule.  So, $\ci(A)$ is a system of  right $W$--modules and $H^*(X;\ci(A))$ ($=H^*_c(\cu)$) is a right $W$--module.  Similarly, $\cac(A)$ is a system of left $W$--modules and $H_*(X;\cac(A))$ ($=H_*(\cu)$) is a left $W$--module. 

For each nonnegative integer $p$, define
\begin{align}
F_p&:=\sum_{|T|\geq p} A^T,	& E_p&:=\bigoplus_{|T|<p} \hA^T, \label{e:F}\\
F'_p&:=\sum_{|T|\geq p} H^T,	& E'_p&:=\bigoplus_{|T|<p} \hH^T, \label{e:F'}
\end{align}
where $|T|:=\Card (T)$.  As in \fullref{s:basics}, these give coefficient systems$\ci(F_p)$ and $\cac(F'_p)$ on $X$.
Note that $F_p$ is a right $W$--module and  $\ci(F_p)$ is a coefficient system of right $W$--submodules of $\ci(A)$.  Similarly,  $\cac(F'_p)$) is a system of left $W$--modules.  (However, $E_p$ and $E'_p$ only have the structure of $\zz$ submodules of $A$.)

\begin{lemma}\label{l:cleandecomp}
We have decompositions (as $\zz$--modules):
\begin{enumerate}
\item
$A=F_p\oplus E_p$ and this induces a decomposition of coefficient systems, $\ci(A)=\ci(F_p)\ \oplus\  \ci(E_p)$.
\item
$A=F'_p\oplus E'_p$ and this induces a decomposition of coefficient systems, 
$\cac(A)=\cac(F'_p)\oplus \cac(E'_p)$.
\end{enumerate}
\end{lemma}

\begin{proof}$\phantom{99}$

(i)\qua By the second formula in \fullref{c:hAT}, $F_p=\bigoplus_{|T|\ge p} \hA^T$; hence, by the first formula in the same corollary, $A=F_p\oplus E_p$.  To get the  decomposition of coefficient systems, we must show that $A^U=(F_p)^U\oplus (E_p)^U$ for all $U\subset S$.  Since $A^U=\bigoplus _{T\supset U} \hA^T$,
	\begin{equation}\label{e:pdecomp}
	A^U=\bigoplus_{\substack{T\supset U\\|T|\ge p}}\hA^T \oplus 
	\bigoplus_{\substack{T\supset U\\|T|< p}}\hA^T .
	\end{equation}
Denote the first summation in \eqref{e:pdecomp} by $B$ and the second one by $C$.
\begin{claim}
$B=(F_p)^U$.
\end{claim}

\begin{proof}[Proof of Claim] 
Obviously, $B\subset (F_p)^U$.  Let $x\in (F_p)^U$.  Since $x\in F_p$, we have
\[
x=\sum_{|T|\ge p} \ga^T,
\]
where $\ga^T\in \hA^T$.  Since $x\in A^U$,
\[
x=\sum_{T\supset U} \gb^T,
\]
where $\gb^T\in \hA^T$.  But $A=\bigoplus_{T\subset S} \hA^T$, so the two decompositions of $x$ coincide.  Therefore, $\ga^T=0$ unless $T\supset U$ and  
\[
x=\sum_{\substack{T\supset U\\|T|\ge p}} \ga^T \in B,
\]
which proves that $(F_p)^U\subset B$.  
\end{proof}
Continuing with the proof of \fullref{l:cleandecomp}, note that
a similar argument shows  $(E_p)^U=C$.  Hence, $A^U=(F_p)^U\oplus (E_p)^U$ and (i) is proved.

(ii)\qua  As before, by \fullref{c:hHT}, $A=F'_p\oplus E'_p$.  To get the decomposition of coefficient systems, we must show that $A_{S-U}=(F'_p)_{S-U}\oplus (E'_p)_{S-U}$ for all $U\subset S$.  Since $A_{S-U}=\bigoplus _{T\subset U} (\hH^T)_{S-U}$,
\begin{equation}\label{e:pdecomp'}
A_{S-U}=\bigoplus_{\substack{T\subset U\\|T|\ge p}}(\hH^T)_{S-U} \oplus 
\bigoplus_{\substack{T\subset U\\|T|< p}}(\hH^T)_{S-U}. 
\end{equation}
Denote the first summation in \eqref{e:pdecomp'} by $B'$ and the second one by $C'$.  We claim that $(F'_p)_{S-U}=B'$.  Obviously, $B'\subset (F'_p)_{S-U}$.   Any $x\in F_p$ can be written in the form
\[
x=\sum_{|T|\ge p} \gamma^T,
\]
where $\gamma^T\in \hH^T$.  Since $\gamma^T\in I_{S-U}$ whenever $T\cap (S-U)\neq \emptyset$,  if $T\not\subset U$, we can set  $\gamma^T=0$ without changing the congruence class of $x$ modulo $I_{S-U}$. So, putting
\[
y= \sum_{\substack{T\subset U\\|T|\ge p}} \gamma^T,
\]
we have $y\equiv x \mod I_{S-U}$ and $y\in B'$.  So, $(F'_p)_{S-U}\subset B'$.  A similar argument shows  $(E'_p)_{S-U}=C'$.  Hence, $A_{S-U}=(F'_p)_{S-U}\oplus (E'_p)_{S-U}$ and (ii) is proved.
\end{proof}

\begin{corollary}\label{c:summand}\hfil
\begin{enumerate}
\item
$F_p\hookrightarrow A$ induces $H^i(X;\ci(F_p))\hookrightarrow H^i(X;\ci(A))$ a $W$--equivariant embedding with image a $\zz$--module direct summand.
\item
$F'_p\hookrightarrow A$ induces $H_i(X;\cac(F'_p))\hookrightarrow H_i(X;\cac(A))$ a $W$--equivariant embedding with image a $\zz$--module direct summand.
\end{enumerate}
\end{corollary}

It follows that $F_{p+1}\hookrightarrow F_p$ induces $H^*(X;\ci(F_{p+1}))\hookrightarrow H^*(X;\ci(F_p))$, an embedding of right $W$--modules.  This gives an associated graded group of right $W$--modules:
\[
H^*(X;\ci(F_p))/H^*(X;\ci(F_{p+1})).
\]
Similarly, we have an embedding
$H_*(X;\cac(F'_{p+1}))\hookrightarrow H_*(X;\cac(F'_p))$ of left $W$--modules and an associated graded group of left $W$--modules.  Our goal in this section is to prove \fullref{t:main} below, which gives  a complete computation of these graded $W$--modules.

For each $T\in \cs$, put 
\[
A^{>T}:=\sum_{U\supset T} A^U\quad\text{and}\quad H^{>T}:=\sum_{U\supset T} H^U.
\]
$A^T/A^\ut$ is a right $W$--module and $H^T/H^\ut$ is a left $W$--module.

\begin{example}\label{ex:sign}
(\textsl{The sign representation})\qua $A^{\emptyset}/A^{>\emptyset}$ is isomorphic to $\zz$ as an abelian group.  We can take the image $\ol{b}_1$ of the basis element $b_1$ ($=e_1$) as the generator.  Since $a_sb_1\in A^{>\emptyset}$, $\ol{b}_1\cdot a_s=0$ for all $s\in S$.  Hence, $\ol{b}_1\cdot s=-\ol{b}_1$.  It follows that $W$ acts on
 $A^{\emptyset}/A^{>\emptyset}$ via the sign representation:
 \[
 \ol{b}_1\cdot  w=(-1)^{l(w)} \ol{b}_1.
 \]
\end{example}

\begin{example}\label{ex:Wfinite}
(\textsl{The case of a finite Coxeter group})\qua  If $W$ is finite, then for any $T\subset S$, $A^T/A^\ut\otimes \qq$ can be identified with a (right) $W$--submodule of the rational group algebra $\qq W$.  Similarly, $H^T/H^\ut\otimes \qq$ is a (left) $W$--submodule of $\qq W$.  L Solomon proved in \cite{sol} that we have direct sum decompositions:
\begin{align*}
\qq W&=\bigoplus _{T\subset S} A^T/A^\ut\otimes \qq \\
\qq W&=\bigoplus _{T\subset S} H^T/H^\ut\otimes \qq
\end{align*}
Of course, there is  no such decomposition over $\zz$.  
For an arbitrary Coxeter group $W$, a similar result for $L^2_\bq (W)$ (the ``$\bq$--weighted $L^2$--completion'' of the regular representation) is proved in \cite[Theorem 9.11]{ddjo}.
\end{example}

\begin{theorem}\label{t:main}
For each nonnegative integer $p$,
\begin{enumerate}
\item
there is an isomorphism of right $W$--modules:
\[
H^*(X;\ci(F_p))/H^*(X;\ci(F_{p+1}))\cong \bigoplus _{|T|=p} H^*(X,X^{S-T})\otimes (A^T/A^\ut);
\]
\item
there is an isomorphism of left $W$--modules:
\[
H_*(X;\cac(F'_p))/H_*(X;\cac(F'_{p+1}))\cong \bigoplus _{|T|=p} H_*(X,X^T)\otimes (H^T/H^\ut).
\]
\end{enumerate}
\end{theorem}

In view of \fullref{cor:Ugpring} and \fullref{l:equisheaf},  this theorem  gives a computation of the $W$--modules associated to the corresponding filtrations of $H^*_c(\cu)$ and $H_*(\cu)$.  To prove the theorem we first need the following lemma.

\begin{lemma}\label{l:psi}
There are isomorphisms of $W$--modules:
\[
\psi\co F_p/F_{p+1}\mapright{\cong} \bigoplus_{|T|=p} A^T/A^\ut \quad \text{and}\quad
\psi'\co F'_p/F'_{p+1}\mapright{\cong} \bigoplus_{|T|=p} H^T/H^\ut.
\]
\end{lemma}

\begin{proof}
The inclusion $A^T\hookrightarrow F_p$ induces a map $A^T\to F_p/F_{p+1}$ and $A^\ut$ is in the kernel; so, we get $A^T/A^\ut\to F_p/F_{p+1}$.  Therefore, we have a map of right $W$--modules:
\[
\phi\co \bigoplus_{|T|=p} A^T/A^\ut\to F_p/F_{p+1}.
\]
By \fullref{c:hAT}, the inclusion $\hA^T\hookrightarrow A^T$ induces an isomorphism (of $\zz$--modules),  $\hA^T\to A^T/A^\ut$.  Also, $F_p = \bigoplus_{|T|=p} \hA^T \oplus F_{p+1}$.  So, we have a commutative diagram (of maps of $\zz$--modules):
\[
\begin{matrix}
\bigoplus_{|T|=p}A^T/A^\ut &&\mapright{\phi}&&F_p/F_{p+1}\cr
&&&&\cr
&\nwarrow&&\nearrow&\cr
&&&&\cr
&&\bigoplus_{|T|=p} \hA^T&&
\end{matrix}
\] 
Since the two slanted arrows are bijections, so is $\phi$.  Therefore,  $\phi$ is an isomorphism of right $W$--modules.  Put $\psi:= \phi^{-1}$.

The definition of the second isomorphism $\psi'$ is similar.
\end{proof}

Here is some more notation.  For any $T\subset S$, put
	\begin{equation*}\label{e:Q}
	Q_{\langle T\rangle}=A^T/A^\ut,\qquad Q'_{\langle T\rangle} = H^T/H^\ut.
	\end{equation*}
Since the right $W$--module $Q_{\langle T\rangle}$ is neither a left $W$--module or even a $\zz$--submodule of a left $W$--module, the definition of its (left) $W_U$--invariants as in \eqref{e:inv} cannot be applied directly.   Similarly, the definition of (right) coinvariants from \eqref{e:coinv} does not apply directly to $Q'_{\langle T\rangle}$.   Nevertheless, for each $U\subset S$,  define:
	\begin{align*}
	(Q_{\langle T\rangle})^U&:= (A^T\cap A^U)/(A^\ut\cap A^U)\\ 	
	(Q'_{\langle T\rangle})_U&:= (H^T)_U/(H^\ut)_U
	\end{align*}
These give coefficient systems of $W$--modules on $X$ defined by
	\begin{align*}
	\ci(Q_{\langle T\rangle})(c)&:=(Q_{\langle T\rangle})^{S(c)}\\
	\cac(Q'_{\langle T\rangle})(c)&:= (Q'_{\langle T\rangle})_{S(c)}
	\end{align*}
respectively.
As in \eqref{e:hAT1} and \eqref{e:hHT1},
	\begin{align}
	(Q_{\langle T\rangle})^U&=
	\begin{cases}
	A^T/A^\ut,	&\text{if $U\subset T$;}\\
	0,	&\text{otherwise}
	\end{cases}\label{e:QTU}\\
	(Q'_{\langle T\rangle})_U&=
	\begin{cases}
	H^T/H^\ut,	&\text{if $U\subset S-T$;}\\
	0,	&\text{otherwise}.
	\end{cases}\label{e:Q'TU}
	\end{align}
	
\begin{lemma}\label{l:QT}\hfil
\begin{enumerate}
\item
$H^i(X;\ci(Q_{\langle T\rangle}))=H^i(X,X^{S-T})\otimes Q_{\langle T\rangle}$
\item
$H_i(X;\cac(Q'_{\langle T\rangle}))=H^i(X,X^T)\otimes Q'_{\langle T\rangle}$
\end{enumerate}
\end{lemma}

\proof$\phantom{99}$

(i)\qua  Using \eqref{e:QTU}, we have
\begin{align*}
C^i(X;\ci (Q_{\langle T\rangle}))&=\{f\co X^{(i)}\to Q_{\langle T\rangle}\mid f(c)=0 \text{ if $c\subset X^{S-T}$}\}\\
&=C^i(X, X^{S-T};Q_{\langle T\rangle})\\
&\cong C^i(X,X^{S-T})\otimes Q_{\langle T\rangle}.
\end{align*}

(ii)\qua Using \eqref{e:Q'TU},
\begin{align*}
C_i(X;\cac(Q'_{\langle T\rangle}))&= \bigoplus_{c\in X^{(i)}} (Q'_{\langle T\rangle})_{S(c)}\\
&=
\begin{cases}
Q'_{\langle T\rangle},		&\text{if $c \not\subset X^T$;}\\
0,			&\text{if $c\subset X^T$;}
\end{cases}\\
&=C_i(X,X^T)\otimes Q'_{\langle T\rangle}.\qquad {\hbox to 0pt{\hglue 0.87in \qed\hss}}
\end{align*}

\begin{lemma}\label{l:exact} For any $U\subset S$, 
\begin{enumerate}
\item  the following sequence of right $W$--modules is exact,
\[
0\mapright{} (F_{p+1})^U\mapright{} (F_p)^U\mapright{\tilde{\psi}}\bigoplus _{|T|=p} (Q_{\langle T\rangle})^U\mapright{} \ 0,
\]
where $\tilde{\psi}$ is the map induced by $\psi$ and
\item  the following sequence of left $W$--modules is exact,
\[
0\mapright{} (F'_{p+1})_{S-U}\mapright{} (F'_p)_{S-U}\mapright{\tilde{\psi}'}\bigoplus _{|T|=p} (Q'_{\langle T\rangle})_{S-U}\mapright{} \ 0,
\]
where $\tilde{\psi}'$ is the map induced by $\psi'$.
\end{enumerate}
\end{lemma}

\begin{proof}
In the proof of \fullref{l:cleandecomp}, in formula \eqref{e:pdecomp}, we showed  
\[
(F_p)^U=\bigoplus_{\substack{|T|\ge p \\ T\supset U}} \hA^T.
\]
Put
\[
B:=\bigoplus_{\substack{|T|= p \\ T\supset U}} \hA^T,
\]
$B$ is a $\zz$--submodule of $(F_p)^U$ and it  maps isomorphically onto $(F_p)^U/(F_{p+1})^U$.  The image of $B$ under $\psi$ is 
\[
\bigoplus_{\substack{|T|= p \\ T\supset U}} A^T/A^\ut =\bigoplus_{\substack{|T|= p \\ T\supset U}} Q_{\langle T\rangle}.
\]
This proves (i).

The proof that the sequence in (ii) is short exact is similar.
\end{proof}

\begin{proof}[Proof of \fullref{t:main}]$\phantom{99}$

(i)\qua By \fullref{l:exact}~(i), we have a short exact sequence of coefficient systems on $X$:
\[
0\mapright{}\ci(F_{p+1})\mapright{}\ci(F_p)\mapright{}\bigoplus_{|T|=p} \ci(Q_{\langle T\rangle})\mapright{}\,0
\]
inducing a short exact sequence of cochain complexes:
\[
0\mapright{}\,C^*(X;\ci(F_{p+1}))\mapright{}\,C^*(X;\ci(F_p))\mapright{}\,\bigoplus_{|T|=p} C^*(X;\ci(Q_{\langle T\rangle}))\mapright{} \ 0.
\]
By the argument for \fullref{c:summand}, $H^*(X;\ci(F_{p+1}))\to H^*(X;\ci(F_p))$ is an injection onto a ($\zz$--module) direct summand.  Hence, the long exact sequence in cohomology decomposes into short exact sequences and we have:
\begin{align*}
H^i(X;\ci(F_p))/H^i(X;\ci(F_{p+1}))&\cong \bigoplus_{|T|=p} H^i(X;\ci(Q_{\langle T\rangle}))\\
&\cong \bigoplus_{|T|=p} H^i(X,X^{S-T})\otimes (A^T/A^\ut)
\end{align*}
where the second isomorphism comes from  
\fullref{l:QT}~(i). 

(ii)\qua  The proof of (ii) is similar.
\end{proof}

\begin{remark}\label{r:spectral}
The decreasing filtration $\supset F_p\supset F_{p+1}\cdots$ of \eqref{e:F} gives a decreasing filtration of cochain complexes 
	\[
	\cdots\supset C^*(X;\ci(F_p))\supset C^*(X;\ci(F_{p+1}))\cdots 
	\]
So, the  quotient cochain complexes have the form 
$C^*(X;\ci(F_p)/\ci(F_{p+1}))$.  Taking homology, we get a spectral sequence with $E_1$--term:
	\[
	E_1^{pq}:=H^{p+q}(X;\ci(F_p)/\ci(F_{p+1}))
	\]
It converges to:
	\begin{equation}\label{e:spectral}
	E_{\infty}^{pq}:=\frac{H^{p+q}(X;\ci(F_p))}{\Ima (H^{p+q}(X;\ci(F_{p+1})))}
	\end{equation}
So, the  import of \fullref{t:main} is  that $E^{pq}_1=E^{pq}_\infty$.
\end{remark}

\section{$H^*(W;\zz W)$}\label{s:H}
Let $K$ denote the geometric realization of the poset $\cs$ of spherical subsets.  (The simplicial complex $K$ is contractible because it is a cone; the cone point corresponds to the minimum element $\emptyset \in \cs$.)  For each $s\in S$, define a subcomplex  $K_s \subset K$ as the geometric realization of $\cs_{\ge \{s\}}$.  Put $\gS:=\cu(W,K)$.  (Alternatively, $\gS$ could be described as the geometric realization of the poset $W\cs$ of all ``spherical cosets,'' ie, the poset of all cosets of the form $wW_T$, with $w\in W$ and $T\in \cs$.)  

By construction $W$ acts properly on $\gS$.  It is proved in \cite{d83} that $\gS$ is contractible. Hence,
\begin{equation*}\label{e:WA}
H^*(W;\zz W)=H^*_c(\gS).
\end{equation*}
As before, $A:=\zz W$.  The filtration $A=F_0\supset\cdots F_p\supset\cdots$ gives $H^*_c(\gS)=H^*(K;\ci(A))$ the structure  of a graded $W$--module.  As in \eqref{e:spectral}, let $E^{pq}_\infty$ is  the right $W$--module in filtration degree $p$ associated to $H^{p+q}(K;\ci(A))$. \fullref{t:main} then has the following corollary.

\begin{theorem}\label{t:main2}
The associated graded group of $H^{p+q}(W;\zz W)$ is given,  as a right $W$--module, by
\[
E^{pq}_\infty=\bigoplus_{|T|=p} H^{p+q}(K,K^{S-T})\otimes (A^T/A^\ut).
\]
\end{theorem}

It follows from \fullref{t:cohomU} that we have a direct sum decomposition of $\zz$--modules:
\[
H^*_c(\gS)\cong \bigoplus_{T\in \cs} H^{*}(K,K^{S-T})\otimes \hA^T.
\]
In view of \fullref{t:main2}, it is natural to conjecture that  $H^*_c(\gS)$  decomposes as above into a direct sum of right $W$--modules. However, in general, there is no such decomposition, as we can see by considering the following example.

\begin{example}\label{ex:free}
Suppose $W$ is the free product of $3$ copies of $\zz/2$.  Then $K$ is the cone on $3$ points.   So, it has $3$ edges.
By \fullref{t:main}, $H^1(K,K^S)\otimes A/A^{>\emptyset}$ is a quotient of $H^1_c(\gS)$.    Let $x\in C^1(K)$ be a cochain ($=$ cocycle) which evaluates to $1$ on one of the edges, call it $c$, and to $0$ on the other two edges.  Choose $s\in S$ which is not a vertex of $c$.  Let $y$ denote the image of $x\otimes 1$ in $H^1(K,K^S)\otimes A/A^{>\emptyset}$.  By \fullref{ex:sign}, $A/A^{>\emptyset}$ has rank $1$ as an abelian group and the $W$--action on it is given by the sign representation.  Hence, $y\cdot s=-y$ in $H^1(K,K^S)\otimes A/A^{>\emptyset}$.  Suppose we had a $W$--equivariant splitting $\gf\co H^1(K,K^S)\otimes A/A^{>\emptyset} \to H^1_c(\gS)$.  When regarded as an element of $C^1_c(\gS)$, $x+x\cdot s$ represents $\gf(y+y\cdot s)$ in $H^1_c(\gS)$, ie, it represents $0$.  But  $x$ and $-x\cdot s$ are not cohomologous cocycles in $C^1_c(\gS)$.  (One can see this by noting that there is a line ($=$ infinite $1$--cycle) on which $x$ evaluates to $1$ and $x\cdot s$ to $0$
; see \fullref{f:nosplit}.)  
Hence, there can be no such splitting $\gf$.
\begin{figure}
\begin{center}
\unitlength=.4mm\small
\begin{picture}(180,100)(0,-20)        
\thinlines      
\put(0,0){\line(1,0){180}}
\put(0,60){\line(1,0){180}}
\put(30,0){\line(0,1){20}}
\put(90,0){\line(0,1){60}}
\put(150,0){\line(0,1){20}}
\put(30,60){\line(0,-1){20}}  
\put(150,60){\line(0,-1){20}}
\multiput(30,0)(30,0){5}{\circle*{2}}
\multiput(30,60)(30,0){5}{\circle*{2}} 
\put(90,30){\circle*{2}}

\linethickness{0.2pt}
\qbezier(90,-22)(75,10)(88,14)
\put(88,14){\vector(2,1){0}}

\qbezier(94,26)(100,30)(94,34)
\put(92,25){\vector(-2,-1){0}}
\put(92,35){\vector(-2,1){0}}

\put(86,-30){${K}$}
\put(98,28){${s}$}
\put(101,2){${c}$}
\put(103,62){${sc}$}

\put(92,-22){\vector(1,2){10}}
\put(88,-22){\vector(-1,2){10}}

\linethickness{1.2pt}
\put(90,0){\line(0,1){30}}
\put(90,0){\line(1,0){30}}
\put(90,0){\line(-1,0){30}}

\end{picture}
\end{center}
\caption{Cocycles $x$ and $x\cdot s$}\label{f:nosplit}
\end{figure}
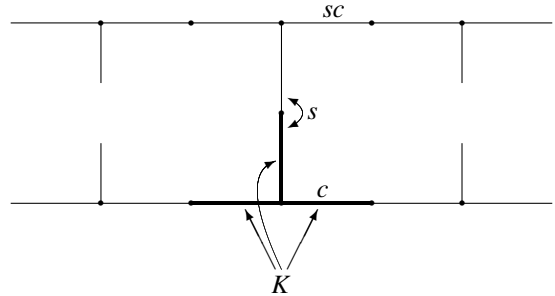
\end{example}

\section{Buildings}\label{s:buildings}

A \emph{chamber system over $S$} is a set $\Phi$ of \emph{chambers} together with a family of equivalence relations on $\Phi$ indexed by $S$.  Two chambers are \emph{$s$--equivalent} if they are equivalent via the equivalence relation with index $s$; they are \emph{$s$--adjacent} if they are $s$--equivalent and not equal.
 A \emph{gallery} in $\Phi$ is a finite sequence of 
chambers $(\gf_0, \dots ,\gf_k)$ such that $\gf_{j-1}$ is adjacent to $\gf_j, 1 \le j \le k$.  The  \emph{type} of this gallery is the word $\bs
=(s_1,  \dots, s_k)$ where $\gf_{j-1}$ is $s_j$--adjacent to $\gf_j$.  If each $s_j$ 
belongs to a given subset $T$ of $S$, then the gallery is a \emph{$T$--gallery}. 
A chamber system is \emph{connected} (resp.\  \emph{$T$--connected}) if any 
two chambers can be joined by a gallery (resp.\ a $T$--gallery).  
The $T$--connected components of a chamber system $\Phi$ are its 
\emph{residues} of  \emph{type $T$}.  Given $\gf\in \Phi$, $\Res(\gf,T)$ denotes the residue of type $T$ containing $\gf$.
 
Suppose $(W,S)$ is a Coxeter system and $M=(m_{st})$  its Coxeter matrix.   
A \emph{building of type $(W,S)$} (or of \emph{type $M$}) is a chamber system $\Phi$ over $S$ such that 
\begin{enumerate}
\item
for all $s\in S$, each $s$--equivalence class contains at least two chambers, and 
\item
there exists a \emph{$W$--valued  distance function}  $\gd\co \Phi \times 
\Phi \to W$.  (This means that given  a reduced word  $\bs$ for an element $w\in W$, chambers $\gf$ and $\gf'$ can be joined by a gallery 
of type $\bs$ if and only if $\gd (\gf, \gf') = w$.   
\end{enumerate}
A residue of type $T$ in a building is itself a building; its type is $(W_T,T)$.  
A building of type $(W,S)$ is \emph{spherical} if $W$ is finite.  A building has \emph{finite thickness} if all $s$--equivalence classes are finite, for all $s\in S$.  (This implies all spherical residues are finite.)
A building is \emph{thick} if all $s$--equivalence classes have at least $3$ elements, for all $s\in S$.

Fix a base chamber $\gf_0\in \Phi$.  The \emph{folding map}, $\pi\co \Phi\to W$, \emph{centered at $\gf_0$} is defined by $\pi(\gf):=\gd(\gf_0,\gf)$.

A Coxeter matrix $M$ is \emph{right-angled} if all its off-diagonal entries are $2$ or $\infty$.  A building is \emph{right-angled} if its Coxeter matrix is.

\begin{example}\label{ex:rafinite}
(\textsl{Right-angled spherical buildings})\qua If $(W,S)$ is right-angled and spherical, then it has the form $W=W_{s_1}\times\cdots\times W_{s_n}$, with $S=\{s_1,\dots s_n\}$ (ie, $W\cong (\zz/2)^S$).
A building of type $(W_{s_i}, \{s_i\})$ is simply a set $\Phi_i$ with at least $2$ elements  (it is thick if it has at least $3$ elements).  It follows that any right-angled spherical building $\Phi$ has the form $\Phi=\Phi_1\times\cdots \times\Phi_n$.  Fix a base chamber $\gf_0=(x_1,\dots,x_n)\in \Phi$, giving us a folding map $\pi\co \Phi\to W$.  Let $L_\Phi$ be the set of $\gf \in \Phi$ such that $\pi(\gf)$ is the longest element in $W$.  Clearly, $L_\Phi=(\Phi_1-x_1)\times\cdots\times (\Phi_n-x_n)$.  So, when $\Phi$ is thick, $L_\Phi$ is also a right-angled spherical building.
\end{example}
In what follows $\Phi$ is a building of finite thickness and  type $(W,S)$. $\pi\co \Phi\to W$ is a folding map.
For the remainder of  this section, $A$ denotes the abelian group of finitely supported ($\zz$--valued) functions on $\Phi$ and   
for any subset $T$ of $S$, $A^T\subset A$ denotes the subgroup of functions which are constant on $T$--residues.  For each $w\in W$, put $\Out(w):=S-\In(w)$, where $\In(w)$ was defined in \fullref{s:gpring}.

\begin{lemma}\label{l:1}
Suppose that $(\Phi,\pi)$ is a spherical building equipped with a folding map. 
For $\gf\in \Phi$, let $g_\gf$ be the characteristic function
of $\Res(\gf,\Out(\pi(\gf)))$. For any $T\subset S$, put $B^T=\{g_\gf\mid \Out(\pi(\gf))\supset T\}$.  Then $B^T$ is a basis of $A^T$.
\end{lemma}

\begin{proof}First we check that $\{g_\gf\}_{\gf\in \Phi}$ is linearly independent.
Suppose we have  a linear relation $\sum_{\gf\in \Phi}\alpha_\gf g_\gf=0$.  Choose $\psi$ with shortest $\pi(\psi)$ 
such that $\alpha_\psi\ne0$. Then $0=(\sum_{\gf\in \Phi}\alpha_\gf g_\gf)(\psi)=\alpha_\psi$, contradiction. 

Put $C^T:=\{\gf\in \Phi\mid \gf \text{ is the shortest element of $\Res(\gf,T)$}\}$.  Then $C^T$ is
a set of representatives for the set of $T$--residues of $\Phi$.  Notice that $T\subset \Out(\pi(\gf))$ if and only if $\gf\in C^T$.
Therefore, $\Card(B^T)=\Card(C^T)=\rank (A^T)$.
Let $e_R$ denote the characteristic function of a residue $R$ in $\Phi$.  Put $R_\gf:=\Res(\gf,T)$. The \emph{standard basis} for $A^T$ is  $\{e_{R_\gf}\}_{\gf\in C^T}$ .    Suppose $T\subset \Out(\pi(\gf))$.  Then $g_\gf=\sum e_R$ where $R$ ranges over the $T$--residues in $\Res(\gf,\Out(\pi(\gf)))$.  In other words, the matrix which expresses $\{g_\gf\}_{\gf\in C^T}$ in terms of the $e_{R_\gf}$ has $1$'s on the diagonal and is upper triangular when $C^T$ is ordered with respect to $l(\pi(\gf))$.  Hence, $B^T$ is also a basis for $A^T$
\end{proof}

Now \emph{suppose $\Phi$ is right-angled, thick and infinite}. Let $R$ be a spherical residue in $\Phi$, and 
$U\subset S$ its type.
Let $L_R$ denote the set of chambers $\gf\in R$ which have the longest possible $\pi(\gf)$.
By \fullref{ex:rafinite}, $L_R$ is a  right-angled spherical building of type $(W_U,U)$. For each spherical residue $R$,  fix a  folding map
$\pi_{L_R}: L_R\to W_U$. Note that $\pi_{L_R}$ is  totally different from (the restriction of) $\pi$. 

\begin{definition}
Let $\gf\in \Phi$, $U=\In(\pi(\gf))$ and $R=R_\gf:=\Res(\gf,U)$. Then $\gf\in L_R$. 
Define $f_\gf$ to be the characteristic function of $\Res_\Phi(\gf,\Out(\pi_{L_R}(\gf)))$,
where $\Out(\ )$ is computed in $W_U$. 
\end{definition}

We claim that the restriction of $f_\gf$ to $L_R$ is the function $g_\gf$ obtained by applying \fullref{l:1}
to $(L_R,\pi_{L_R})$. 
To see this, it is enough to check the following fact: 
for any $T\subset U$, 
the intersection of a $T$--residue in $\Phi$ with $L_R$ is either
empty or is a $T$--residue in $L_R$. 
Let $\psi\in L_R$ and let
$p\co R\to W_U$ be the $\psi$--based folding map. Then
$p$ maps $\Res_R(\psi,T)$ and $\Res_{L_R}(\psi,T)$ onto $W_T$, while
any other $T$--residue in $L_R$ is mapped onto a nontrivial
coset of $W_T$ in $W_U$. It follows that $\Res_R(\psi,T)\cap L_R=\Res_{L_R}(\psi,T)$.

Observe that if $g_\gf$ is constant on $T$--residues in $L_R$, then $T\subset \Out(\pi_{L_R}(\gf))$ and hence,  $f_\gf$ is constant
on $T$--residues in $\Phi$. Also, the support of $f_\gf$ is contained in $R_\gf$.

Not every spherical residue is of the form $R_\gf$. Call a spherical residue $R$  \emph{greedy},
if for some (hence, every) $\gf\in L_R$, we have $R=\Res(\gf,\In(\pi(\gf)))$. For a greedy residue
$R$ of type $U$ and a subset $T\subset U$, put $B^{R,T}=\{f_\gf\mid \gf\in L_R,\  \Out(\pi_{L_R}(\gf))\supset T\}$
(for convenience, put $B^{R,T}=\emptyset$ if $T\not\subset U$).
Then $\{f|_{L_R}\mid f\in B^{R,T}\}$ coincides with the basis $B^T$ associated 
to $(L_R,\pi_{L_R})$  in \fullref{l:1}. Therefore, $B^{R,T}$ is linearly independent.
In particular, the set $B^{R,\emptyset}=\{f_\gf\mid \gf\in L_R\}$ is linearly independent.

For any spherical $T\subset S$ we put $B^T=\bigcup B^{R,T}$, where the union is over all
greedy residues. Also put $B=B^\emptyset=\{f_\gf\mid \gf\in \Phi\}$.  The next result is the analog of \fullref{l:bbasis}.

\begin{theorem}\label{t:basisbldg}
For each $T\in\cs$, $B^T$ is a basis of $A^T$.
\end{theorem}

\begin{proof}
First we check that $B$ is linearly independent.
Let $\sum_{\gf\in \Phi}\alpha_\gf f_\gf=0$.
Suppose $\psi\in \Phi$ is an element with $\alpha_\psi\ne0$ and maximum $l(\pi(\psi))$.
Let $R=\Res(\psi,\In(\pi(\psi)))$. Then $\sum_{\gf\in \Phi}\alpha_\gf f_\gf|_{L_R}=
\sum_{\gf\in L_R}\alpha_\gf f_\gf|_{L_R}$. Since the set $\{f_\gf\mid \gf\in L_R\}$ is a linearly independent,  we have $\alpha_\gf=0$ for $\gf\in L_R$.
In particular $\alpha_\psi=0$, contradiction.

Suppose now that $f\in A^T$ is not a linear combination of elements of $B^T$. Choose such an $f$
with minimum $M:=\max\{l(\pi(\gf))\mid f(\gf)\ne 0\}$. Let $\psi\in \Phi$ be such that 
$f(\psi)\ne0$ and $\ell(\pi(\psi))=M$. Let $R=\Res(\psi,\In(\pi(\psi)))$.
Clearly, $U=\In(\pi(\psi))\supset T$. 
Let $(R_i)_{i=0,\ldots,N}$ be the collection of $U$--residues in $\Phi$ that are $U$--connected
components of $\pi^{-1}(\pi(R))$ (where $R_0=R$). These residues are pairwise disjoint.  We also note that they are greedy.  Indeed, (1) $R_0$ is greedy by its definition, (2) the notion of greediness for a residue depends only on its image under $\pi$ and (3) $\pi(R_i)=\pi (R_0)$.  
The restrictions to $L_{R_i}$ of functions
from $B^{R_i,T}$ form a basis of the space of functions on $L_{R_i}$ constant
on $T$--residues -- in fact, this is exactly the basis produced by applying \fullref{l:1} to $L_{R_i}$.
Since the restriction of $f$ to $L_{R_i}$ is constant on $T$--residues, there exists $f_{R_i}\in \Span(B^{R_i,T})$
whose restriction to $L_{R_i}$ coincides with that of $f$.
Let $\widetilde{f}=f-\sum_{i=0}^Nf_{R_i}$. Then $\widetilde{f}$
is in $A^T$, is not in $\Span(B^T)$ and has value $0$ on $\pi^{-1}(\pi(\psi))$. 
Observe also that if $\gf\not\in\pi^{-1}(\pi(\psi))$ and $\ell(\pi(\gf))=M$, then 
$\widetilde{f}(\gf)=f(\gf)$. Therefore, we can repeat our procedure for all such $\gf$, and
finally obtain a counterexample to $A^T=\Span(B^T)$ with smaller $M$, a contradiction. 
\end{proof}  

Put $\wh{B}^T:=B^T-\bigcup_{U\supset T} B^U$ and $\widehat{A}^T:=\Span(\wh{B}^T)$.  Just as in \fullref{s:gpring}, \fullref{t:basisbldg} has the following corollary.  
\begin{corollary}\label{c:ATbldg}
\textup{(Compare \fullref{c:hAT})} 
\[
A^U=\bigoplus_{T\supset U}\widehat{A}^T.
\]
\end{corollary}

As in \fullref{ex:simple},  given a mirror structure $(X_s)_{s \in S}$ on a CW complex $X$, we get  a covariant coefficient system $\ci(A)$ on $X$ defined by $\ci(A)(c):=A^{S(c)}$. As in \eqref{e:Asheaf}, \fullref{c:ATbldg} means that we have a decomposition of coefficient systems:
\begin{equation}\label{e:decompbldg}
\ci(A)=\bigoplus_{T\in \cs} \ci(\wh{A}^T),
\end{equation}
where $\ci(\wh{A}^T)(c):=\wh{A}^T\cap A^{S(c)}$.

\medskip
\textbf{The geometric realization of a building}\qua
Given a CW complex $X$ with mirror structure,
define an equivalence relation  $\sim$ on $\Phi \times X$ by $(\gf, x) \sim (\gf', x')$ if and 
only if $x = x'$ and $\gd (\gf, \gf') \in W_{S(x)}$.  The $X$--\emph{realization} 
of $\Phi$, denoted  $\cu(\Phi,X)$, is defined by
	\begin{equation}\label{e: defUPhi}
	\cu(\Phi,X):=(\Phi \times X)/\sim.  
	\end{equation}
($\Phi$ has the discrete topology.)  When $X=K$ (the geometric realization of $\cs$), $\cu(\Phi, K)$ is denoted by $|\Phi|$ and called the \emph{geometric realization of $\Phi$}. 

As in \fullref{cor:Ugpring} and \fullref{l:equisheaf}, identify $C^*(\cu(\Phi,X))$ with $C^*(X;\ci(A))$.  
Using \eqref{e:decompbldg}, the proof of \fullref{t:cohomU} goes through to give the following.

\begin{theorem}\label{t:cohombldg}
Suppose $\Phi$ is a right-angled, thick, infinite building. Then
\[
H^i_c(\cu(\Phi,X))\cong \bigoplus_{T\in \cs} H^i(X,X^{S-T})\otimes \hA^T.
\]
In particular,
\[
H^i_c(|\Phi|)\cong \bigoplus_{T\in \cs} H_i(K,K^{S-T})\otimes \hA^T.
\]
\end{theorem}

\begin{Remark}
A similar result for any building of finite thickness is claimed in \cite[Theorem 5.8]{dmei}; however, there is a mistake in the proof.
\end{Remark}

\section{Hecke algebra coefficients}\label{s:hecke}
In this section we work over the rational numbers $\qq$ rather than $\zz$.

Let $i:S\to I$ be a function to some index set $I$ such that $i(s)=i(s')$ whenever $s$ and $s'$ are conjugate in $W$.  Let $\bq=(q_i)_{i\in I}$ be a fixed $I$--tuple of 
rational numbers.  Write $q_s$ instead of $q_{i(s)}$.  If $s_1\cdots s_l$ is a reduced expression for an element $w\in W$, then the number  $q_{s_1}\cdots q_{s_l}$ is independent of the choice of reduced expression.  We write it as $q_w$.  The \emph{Hecke algebra} $A_\bq$ of $W$
is a deformation of  the group algebra $\qq W$ which is equal to $\qq W$ when each $q_s=1$.  As a rational vector space, it has the same basis $\{e_w\}_{w\in W}$ as does $\qq W$.  Multiplication is determined by the rules:
	\begin{align*}
	e_we_{w'}&=e_{ww'},\quad \text{if $l(ww')=l(w)+l(w')$}\\
	e_s^2&=(q_s-1)e_s+q_s.
	\end{align*}
	
Given a special subgroup $W_T$, $A_\bq(W_T)$ denotes the Hecke algebra of $W_T$.  It is a subalgebra of $A_\bq$.  There are ring homomorphisms $\ga\co A_\bq(W_T)\to \qq$ and $\gb\co A_\bq(W_T)\to \qq$, defined by $\ga(e_w):=q_w$ and $\gb(e_w):=(-1)^{l(w)}$, respectively.  Given a left $A_\bq$--module $M$ and a subset $T$ of $S$, put 
\[
M^T:=\{x\in M\mid ax=\ga(a)x \text{ for all } a\in A_\bq(W_T)\}.
\]
This gives a coefficient system $\ci(M)$ on $X$ in the same way as  Examples~\ref{ex:simple}.

As in \cite{ddjo}, for each $T\in \cs$, we modify the formulas in \eqref{e:aThT} to define elements $a_T$ and $h_T$ in $A_\bq$ by
	\begin{equation*}\label{e:aThT2}
	a_T:=\frac{1}{W_T(\bq)}\sum_{w\in W_T} e_w\quad \text{and} 
	\quad h_T:=\frac{1}{W_T(\bq^{-1})}\sum_{w\in W_T} (-1)^{l(w)}q_w^{-1}e_w
	\end{equation*}
where
\[
W_T(\bq):=\sum_{w\in W_T} q_w,\quad \text{and}\quad 
W_T(\bq^{-1}):=\sum_{w\in W_T} q_w^{-1}.
\]
Put $A_\bq^T:=a_TA_\bq$, $H_\bq^T:=A_\bq h_T$.  (If $T\notin \cs$,  $A_\bq^T:=0$, $H_\bq^T:=0$.)

For each subset $U$ of $S$, put
\[
(A_\bq)_U:= A_\bq\otimes_{A_\bq(W_U)} \qq=A_\bq/A_\bq I_U,
\]
where $A_\bq (W_U)$ acts on $\qq$ via the symmetric character $\ga_U$ and $I_U:=\sum_{s\in U} H^s_\bq$ is the augmentation ideal of $A_\bq(W_U)$.  $A_\bq^U$ is a right $A_\bq$--module and $(A_\bq)_U$ is a left $A_\bq$--module.

We have decreasing filtrations $(F_p)$ and $(F'_p)$ of $A_\bq$, defined exactly as in \eqref{e:F} and \eqref{e:F'}.

If  $X$ and $\cu$ are as before, then the proof of \fullref{t:main} gives the following.

\begin{theorem}\label{t:main3}
With notation as above, for each nonnegative integer $p$,
\begin{enumerate}
\item
there is an isomorphism of right $A_\bq$--modules:
\[
H^*(X;\ci(F_p))/H^*(X;\ci(F_{p+1}))\cong \bigoplus _{|T|=p} H^*(X,X^{S-T})\otimes (A^T_\bq/A^\ut_\bq);
\]
\item
there is an isomorphism of left $A_\bq$--modules:
\[
H_*(X;\cac(F'_p))/H_*(X;\cac(F'_{p+1}))\cong \bigoplus _{|T|=p} H_*(X,X^T)\otimes (H^T_\bq/H^\ut_\bq).
\]
\end{enumerate}
\end{theorem}

\noindent
\textbf{$\boldsymbol{BN}$ pairs}\qua The importance of Hecke algebras
lies in their relationship to buildings and $BN$ pairs (eg, see
Bourbaki \cite[Exercises 22 and 24, pages 56--58]{bourbaki}).  Suppose
that $(G,B)$ is a $BN$ pair.  Associated to $(G,B)$ we have a Coxeter
system $(W,S)$ and for each $s\in S$, a subgroup $G_s$ of $G$ such
that $G_s = B \cup BsB$.  Put $\Phi:=G/B$.  For each subset $T$ of
$S$, put $G_T:=BW_T B$.  Two cosets $gB$ and $g'B$ are
$s$--equivalent if they have the same image in $G/G_s$.  This gives
$\Phi$ the structure of a building.  $\Phi$ has finite thickness if
$(G_s:B)<\infty$ for all $s\in S$.  If this is the case, put
$q_s=(G_s:B)-1$ and regard $\bq=(q_s)$ as an $I$--tuple, where $I$ is
the set of conjugacy classes of elements in $S$.

Let $F(G/B)$ denote the $\qq$ vector space of finitely supported, $\qq$--valued functions on $G/B$.  The left $G$--action on $G/B$ gives $F(G/B)$ the structure of a right $G$--module.  For any $T\subset S$, we have a projection $p_T\co G/B\to G/G_T$.  By pulling back via $p_T$, we can identify $F(G/G_T)$ with a $G$--submodule of $F(G/B)$.

Regard $F(G/B)$ as a subset of all $\qq$--valued functions on $G$. 
The Hecke algebra can be identified with the subspace of $F(G/B)$ consisting of those functions which are invariant under the $B$--action on $F(G/B)$ (induced from the left $B$--action on $G/B$).  Under this identification, the basis element $e_w\in A_\bq$ is identified with the characteristic function of the double coset, $C(w^{-1}):=Bw^{-1}B$ and the idempotent $a_T$ with the characteristic function of $G_T$ ($=BW_TB$).

Given $f\in F(G/B)$ and  $a\in A_\bq$ their \emph{convolution} is defined by
	\[
	(a*f)(h):=\int_G a(g^{-1}h)f(g)dg.
	\]
Here we are integrating with respect to Haar measure $dg$ normalized so that the measure of $B$ is $1$.  So, $A_\bq$ acts from  the left on $F(G/B)$ by convolution,  (In fact, $A_\bq$ is the intertwining algebra ($=$ the commutant) of $G$ on $F(G/B)$.)

\begin{lemma}\label{l:tensor}
For any $T\in \cs$,
\[
A^T_\bq\otimes_{A_\bq} F(G/B)=F(G/G_T).
\]
\end{lemma}
\begin{proof}
Let $\sum_{i} a_T\ga_i\otimes f_i$ be a typical element of the left-hand side.  It can be rewritten as $1\otimes \sum a_T*\ga_i *f_i $.  The universal map to $F(G/B)$ consists of taking the second factor.  Since $a_T$ is the characteristic function of $G_T$, $a_T*\ga_i*f_i$ lies in   $F(G/G_T)$; so, the image of this map is the right-hand side.
\end{proof}

\begin{theorem}\label{t:bldgcochains}
\[
C^i_c(|\Phi|;\qq)=C^i(K;\ci(A_\bq))\otimes_{A_\bq} F(G/B)
\]
\end{theorem}

\begin{proof}
This is a direct consequence of \fullref{l:tensor}, since
$$
C_c^i(|\Phi|;\qq)=\bigoplus_{c\in K^{(i)}} F(G/G_{S(c)}).
\proved$$
\end{proof}

The natural conjecture is the following.
\begin{conjecture}\label{conj:7.3}
\[
H^i_c(|\Phi|)=H^i(K;\ci(A_\bq))\otimes_{A_\bq} F(G/B).
\]
\end{conjecture}

The filtration $(F_p)$ induces a filtration of $H^*(K;\ci(A_\bq))$ and hence, of $H^*_c(|\Phi|)$.  So,  \fullref{l:tensor} leads us to the following.
\begin{conjecture}\label{conj:7.4}
In filtration degree $p$, the associated graded group of $H^*_c(|\Phi|)$ is given, as a right $G$--module, by
\[
\bigoplus_{|T|=p} H^*(K,K^{S-T})\otimes F^T/F^\ut,
\]
where $F^T:=F(G/G_T)$ and $F^\ut$ denotes the submodule spanned by the $F^U$ with $U\supsetneq T$.
\end{conjecture}

\bibliographystyle{gtart}
\bibliography{link}

\end{document}